\documentclass[12pt]{article}

\usepackage{amsmath}
\usepackage{amsthm}
\usepackage{graphicx}
\usepackage{amssymb}
\usepackage{float}
\usepackage{subcaption}
\usepackage{dsfont}
\usepackage{url}
\usepackage{epsfig}
\usepackage{epsf}
\usepackage{multicol}
\usepackage{titletoc}
\usepackage{geometry}
\usepackage{dsfont}
\usepackage{multirow}
\usepackage{sectsty}
\usepackage{array}
\usepackage{booktabs}

\newcommand{\beqn}{\begin{eqnarray*}}
\newcommand{\eeqn}{\end{eqnarray*}}
\newcommand{\bneqn}{\begin{eqnarray}}
\newcommand{\eneqn}{\end{eqnarray}}
\newcommand{\parens}[1]{\left(#1\right)}

\newcommand{\bracks}[1]{\left[#1\right]}
\newcommand{\abss}[1]{\left|#1\right|}

\renewcommand{\exp}[1]{\mathrm{exp}\parens{#1}}

\newcommand{\convp}{~{\buildrel p \over \rightarrow}~}

\newcommand{\convd}{~{\buildrel \mathcal{D} \over \rightarrow}~}
\def\Cox{\hfill \Box}

\newtheorem{remark}{Remark}[section]
\newtheorem{theorem}{Theorem}[section]
\newtheorem{proposition}{Proposition}[section]
\newtheorem{corollary}{Corollary}[section]
\newtheorem{lemma}{Lemma}[section]
\newtheorem{definition}{Definition}[section]

\title{Asymptotics for the time of ruin in the war of attrition}
\author{Philip A. Ernst \footnote{Department of Statistics, Rice University} and Ilie Grigorescu \footnote{Department of Mathematics, University of Miami}}

\begin{document}
\maketitle

\begin{abstract}
We consider two players, starting with $m$ and $n$ units, respectively. In each round, the winner is decided with probability proportional to each player's fortune, and the opponent loses one unit. We prove an explicit formula for the probability $p(m,n)$
that the first player wins.  
When
$m\sim Nx_{0}$, $n\sim N y_{0}$, 
we prove the fluid limit as $N\to \infty$. When
$x_{0}=y_{0}$, then 
$z\to p(N,N+z\sqrt{N})$ converges to the standard normal cumulative distribution function (CDF) and the difference in fortunes scales diffusively. The exact limit of the time of ruin $\tau_{N}$ 
is established as $(T-\tau_N) \sim N^{-\beta}W^{\frac{1}{\beta}}$, $\beta=\frac{1}{4}$, $T=x_{0}+y_{0}$. Modulo a constant, $W \sim \chi^{2}_{1}(z_{0}^{2}/T^{2})$. \\
\end{abstract}


\section{Introduction}
In this paper, we develop a ruin problem which derives inspiration from two important branches of literature in the applied probability community and which was brought to our attention by Robert W. Chen and Larry Shepp. The first of the two is the renowned ``War of Attrition,'' a game theoretic model due to John Maynard Smith (\cite{Smith2}, \cite{Smith1}), and later generalized by \cite{Bishop}, and which has proven to be critical for understanding animal conflict and behavior, particularly within the context of evolutionary stable strategies (see \cite{Haigh}, \cite{Hines}, and \cite{Taylor}, among others). The second of the two is the classical gambler's ruin problem, dating to the work of de Moivre in 1711 (\cite{DeMoivre}), and which has been extensively studied and furthered in the works of \cite{Dubins}, \cite{Feller}, \cite{Hald}, \cite{Ross}, \cite{Uspensky}, and \cite{Whitworth}, among others. A further modification of the classical gambler's ruin recently appeared in \cite{Katriel}. \\
\indent In all of the above two-player scenarios, the losing player must give one unit (a point, a dollar, etc.) to his opponent. A very compelling reformulation of this appeared in 1979 in \cite{Kaigh}, and the author called it the ``the attrition ruin problem.'' In the setup of \cite{Kaigh}, the two opponents are in attrition and aim to wear the other down over time. However, unlike the classical gambler's ruin, when the losing player loses a point, the point does not go to his opponent; the point is simply discarded and the winning player stays as is. The author of \cite{Kaigh} claims that this model is in many ways better suited for modeling games between opponents in contests such as board games and in best-of-seven series in sports. Yet, despite this model's apparent applicability, it has only been mentioned once in the literature (\cite{Kozek1}). We surmise that the reason is that the work of (\cite{Kaigh}) is limited by one of its key assumptions: there is a constant probability of winning and losing on each turn.\\
\indent We thus formulate the following model below, which, like the model in \cite{Smith2}, involves a ``war'' between units, and like the model in \cite{Kaigh}, is an attrition ruin problem. We call it the ``war of ruins.'' 
The game is played by two armies, which for convenience we call army A and army B, starting with a pair $(m,n)$ of
soldiers, where $m$ and $n$ are non-negative integers, designating the assets of armies A and B respectively. At discrete times, as long as $m, n>0$, army A (with $m$ units) wins with probability $p=m/(m+n)$ and its number of units remains the same. When this happens, Army B loses one unit so that $(m,n)\to(m, n-1)$. Further, with probability $n/(m+n)$, army B wins and $(m,n)\to(m-1, n)$. Let $p(m,n)$ the ``ruin'' probability; namely, the probability that army A is reduced to 0 soldiers before Army B is reduced to 0 soldiers.  \\
\indent It is important to note that our model is in some ways similar to the work on ruin probabilities for correlated insurance claims.
For references to this very active  literature, we refer to \cite{Alb}, \cite{Asmussen1}, \cite{Asmussen2}, and \cite{Mikosch}, among others.\\
\indent In Section \ref{sec2}, we prove an explicit formula for $p(m,n)$:

\bneqn
p(m,n)=\sum_{j=0}^n \frac{(-1)^j}{j!} \frac{(n-j)^{m+n}}{(m+n-j)!}, \,\, \,\, m\geq 0, n\geq 0, m+n >0.
\eneqn
We then present our first key result (Theorem \ref{thm1}): as $m \rightarrow \infty$,
\bneqn \label{neceq1}
p(m,m+x\sqrt{m})\rightarrow \Phi(x),\,\, -\infty<x<\infty,
\eneqn
where $\Phi$ is the standard normal cumulative distribution function (CDF). A difficult and essential part of the theorem is completed by noticing a surprising connection between $p(m,n)$ and the Eulerian numbers. In Section \ref{sec3}, we take a closer look at the model of \cite{Kaigh}, where $A$ and $B$ lose a unit with fixed probabilities independent of $m$ and $n$. The models studied can be viewed as random walks in the first quadrant, for which we refer the reader to \cite{Kill2}, \cite{Kill3}, and \cite{Kill1}. In Section \ref{sec4}, we consider the continuous time version of the model of Section \ref{sec2} with $X_{t}$, $Y_{t}$ to be the assets of the two players at time $t\ge 0$ and ${\mathcal T}(m,n)$ be the duration of the game, i.e. the time when 
one of the players assets are reduced to zero, assuming exponentially distributed interevent times.  We prove a {\em fluid limit} in Theorem \ref{t_det} based on the law of large numbers scale. Namely, if the two players start with $m=X^{N}_{0} \sim N x_{0} $, $n=Y^{N}_{0}\sim N y_{0}$, $T=x_{0}+y_{0}$ and $t\to Nt$, then  $N^{-1}(X^{N}_{Nt}, Y^{N}_{Nt})$ converges in distribution, as $N\to \infty$, to a pair of coupled solutions of time inhomogeneous differential equations exploding in a finite (nonrandom) time $\tau \le T$.
The result reveals that $x_{0}-y_{0}=0$ is critical, implying that $\tau=T$ and
a finer scaling (diffusive) is available, 
when $Z^{N}_{0}=X^{N}_{0}-Y^{N}_{0}=N^{\frac{1}{2}} z_{0}$. The difference scales 
under $t\to N t$ (Theorem \ref{thmd})
to a diffusion bearing similarities to the Brownian bridge.
If $\tau_{N}=N^{-1}{\mathcal T}(X^{N}_{0},Y^{N}_{0})$ is the scaled time of ruin, then
Theorem \ref{tsc} (the second key result of our work)
and Corollary \ref{csc} determine exactly the limiting distribution of the 
residual time $T-\tau_{N} \sim N^{-\beta} W^{\frac{1}{\beta}}$, where (modulo a known constant)
$W\sim \chi^{2}_{1}(\frac{3z^{2}}{T})$, the non-central $\chi^{2}$ with one degree of freedom and non-centrality parameter 
$3z^{2}/T$. When $z=0$ this is simply $\chi^{2}_{1}$. Remarkably, if $T-\tau_{N} $ is seen as a fluctuation term from the deterministic limit $T$, then the scale is non-Gaussian (which would correspond to $\beta=\frac{1}{2}$ as in the classical CLT), being equal to $\beta=\frac{1}{4}$ instead.

The idea of the proof is to determine a sufficiently large 
family of martingales for the limiting diffusion, in order to evaluate the moments of the residual time. These are obtained as confluent hypergeometric functions 
indexed by a continuous parameter in (\ref{k2}). The most challenging technical difficulty, as is the case with scaling limits, is the replacement of the smooth martingales with their discrete approximations. This is done in Section \ref{s:tsc}, where 
careful estimates near the singularity $\tau_{N}$ are carried out.

\section{Ruin probabilities} \label{sec2}
For $m>0$ and $n>0$, the recurrence
\bneqn \label{eq1}
p(m,n)=\frac{n}{m+n}p(m-1,n)+\frac{m}{m+n}p(m,n-1)
\eneqn
immediately holds.
The corresponding boundary conditions are
\bneqn \label{eq2}
p(m,0)=0, \,\, m>0\,\, ;\,\, p(0,n)=1, \, n>0.
\eneqn
We now find an explicit form of $p(m,n)$ in Proposition \ref{prop1}.
\begin{proposition} \label{prop1}
\bneqn \label{eq3}
p(m,n)=\sum_{j=0}^n \frac{(-1)^j}{j!} \frac{(n-j)^{m+n}}{(m+n-j)!}\,\,, \,\, \,\, m\geq 0, n\geq 0,\, m+n >0,
\eneqn
satisfies (\ref{eq1}) and (\ref{eq2}) which uniquely determine $p(m,n)$.
\end{proposition}
\begin{proof}
Proof by induction. First, note that (\ref{eq3}) agrees with (\ref{eq2}) for $m=0$. In addition, (\ref{eq3}) agrees with (\ref{eq2}) for $n=0$. This completes the base case. We now proceed with the induction step. As it is defined in (\ref{eq3}), $p(m,n)$ must satisfy (\ref{eq1}) since
\bneqn \label{eq222}
& &\frac{n}{m+n}p(m-1,n)+\frac{m}{m+n}p(m,n-1) \nonumber \\
&=&\frac{1}{(m+n)!}\sum_{j=0}^{n-1} {m+n-1 \choose j} (-1)^j \bracks{n(n-j)^{m+n-1}+m(n-j-1)^{m+n-1}} \nonumber \\
&=& \frac{1}{(m+n)!} \sum_{j=0}^{n-1}\bracks{{m+n-1 \choose j}n-{m+n-1 \choose j-1}m}(-1)^j (n-j)^{m+n-1} \nonumber \\
&=& \sum_{j=0}^n \frac{(-1)^j}{j!}\frac{(n-j)^{m+n}}{(m+n-j)!},
\eneqn
where we have used the combinatorial identity
\beqn
{m+n-1 \choose j}n- {m+n-1 \choose j-1}m=(n-j){m+n \choose j}.
\eeqn
Since (\ref{eq222}) agrees with $p(m,n)$ in (\ref{eq3}), (\ref{eq1}) follows. This completes the proof. $\Cox$
\end{proof}
Note that equation (\ref{eq3}) is not suitable for analysis with large $n$. We thus turn to examining $p(m,n)$ for $m$ and $n$ both large.  We begin with Lemma \ref{lemer1} below. As we shall see, the closed form of equation (\ref{thmeq1}) in Lemma \ref{lemer1} is essential for proving the asymptotic formula (\ref{neceq}) in Theorem \ref{thm1}.
\begin{lemma} \label{lemer1}
For $(x,y) \in [0,1) \times [0,1),$ 

\bneqn \label{thmeq1}
\phi(x,y)= \sum_{m=0}^\infty \sum_{n=0}^\infty p(m,n)x^m y^n= \frac{xe^{-x}}{xe^{-x}-ye^{-y}}+\frac{y}{1-y}\frac{1}{y-x}
\eneqn
holds for the generating function $\phi(x,y)$.
\end{lemma}
\noindent Before supplying the proof, we pause to mention a key idea: we``guessed'' the correct right-hand side of the above equation by noticing, for small $m$ and $n$, that $p(m,n)$ are related to the Eulerian numbers $A_{m,n}$ (see \cite{Stanley}) by
\bneqn \label{inter}
A_{n+m,n}=(n+m)!\parens{p(m,n)-p(m+1,n-1)},\,\, m>0.
\eneqn
\begin{proof}
To calculate $p(m,n)$ for $m$ and $n$ both large, consider the generating function appearing after the first equality of display (\ref{thmeq1}) above.
Without loss of generality, $p(0,0)=1$.
After lengthy algebraic manipulation, it can be shown, using equation (\ref{eq1}), that
\bneqn \label{eq6}
x(1-y)\phi_x+ y(1-x)\phi_y=\frac{y}{(1-y)^2},\,\, \phi_x=\frac{\partial{\phi}}{\partial{x}},\,\,\phi_y=\frac{\partial{\phi}}{\partial{y}}
\eneqn
holds.
Further lengthy algebraic manipulation can now be employed to show that
\bneqn \label{eq7}
\psi(x,y)=\frac{xe^{-x}}{xe^{-x}-ye^{-y}}+\frac{y}{1-y}\frac{1}{y-x},
\eneqn
also satisfies equation (\ref{eq6}) when $\phi$ is replaced by $\psi$. We now let $(x_0,y_0)$ be an arbitrary point of $(0,1) \times (0,1)$ and let, for $t \leq 0$, $x=x_t$, $y=y_t$ be trajectories defined by
\bneqn \label{eq8}
\dot{x}_t=x(1-y),\,\, \dot{y}_t=y(1-x).
\eneqn
The derivative of $\phi(x_t,y_t)$-$\psi(x_t,y_t)$ is zero along the trajectory, as $\phi$ and $\psi$ satisfy (\ref{eq6}).
Along the trajectories of (\ref{eq8}), for some constant $c>0$, 
\bneqn \label{eq9}
x_te^{-x_t}=cy_te^{-y_t}
\eneqn
must hold.
Since $x_0(1-y_0)>0$ and $y_0(1-x_0)>0$, as $t$ decreases below zero, then $x_t$ and $y_t$ will also decrease. Thus, for some $t_0<0$, $x_{t_0}=y_{t_0}=0$. Given that $\phi(x_t, y_t)-\psi(x_t,y_t)$ is constant along the trajectory, we have for any arbitrary point $(x_0,y_0)$ that
\bneqn \label{eq10}
\quad \phi(x_0,y_0)-\psi(x_0,y_0)=\phi(x_{t_0},y_{t_0})-\psi(x_{t_0}, y_{t_0})= \phi(0,0)-\psi(0,0)=0,
\eneqn
where $\phi(0,0)=\psi(0,0)=0$. Thus, $\phi \equiv \psi$, and for $(x,y) \in [0,1) \times [0,1),$ equation (\ref{thmeq1}) holds. This concludes the proof. $\Cox$
\end{proof}

We now proceed to prove our first key result (Theorem \ref{thm1}).

\begin{theorem} \label{thm1} As $m \rightarrow \infty$,
\bneqn \label{neceq}
p(m,m+x\sqrt{m})\rightarrow \Phi(x),\,\, -\infty<x<\infty,
\eneqn
where $\Phi$ is the standard normal CDF.
\end{theorem}
\begin{proof}
It seems intractable to proceed directly from (\ref{eq3}) since the terms become large in absolute value. In lieu, we proceed with L\'{e}vy's convergence theorem. To utilize L\'{e}vy's convergence theorem, we first define the following characteristic functions:
\bneqn \label{eq11}
\phi_m(u)=\sum_{n=0}^\infty \bracks{p(m,n)-p(m,n-1)}e^{inu}= E(e^{i\xi_m u}),
\eneqn
where the $\xi_m$ are random variables such that
\beqn
P(\xi_m \leq n)=p(m,n).
\eeqn
To prove the theorem, we must show that the distribution of $\frac{(\xi_m-m)}{\sqrt{m}}$
is asymptotically normal. Equivalently, by L\'{e}vy's convergence theorem,  we must show that for $-\infty<\theta<\infty$,
\bneqn \label{eq15}
\phi_m\parens{\frac{\theta}{\sqrt{m}}}e^{-i\theta\sqrt{m}}= E(e^{\frac{i\theta\parens{\xi_m-m}}{\sqrt{m}}})\rightarrow e^{-\frac{\theta^2}{2}}.
\eneqn
From equation (\ref{thmeq1}), we observe that:
\bneqn \label{eq12}
\phi(x,y)(1-y)&=&\sum_{m=0}^\infty \sum_{n=0}^\infty \bracks{p(m,n)-p(m,n-1)}x^my^n \nonumber \\
&=& \sum_{m=0}^\infty \phi_m(u)x^m, \,\, y=e^{iu} \nonumber \\
&=&\frac{xe^{-x}(1-y)}{xe^{-x}-ye^{-y}}+\frac{y}{y-x}.
\eneqn
Letting $x=z$, we employ Cauchy's formula in equation (\ref{eq12}) and obtain
\bneqn \label{eq13}
\phi_m(u)=\frac{1}{2\pi i}\int_C \frac{1}{z^{m+1}}\bracks{\frac{ze^{-z}(1-e^{iu})}{ze^{-z}-e^{iu}e^{-e^{iu}}}+\frac{e^{iu}}{e^{iu}-z}}dz,
\eneqn
where $C$ is a contour surrounding zero, with $\abss{z}<1$ on $C$. We now shift the contour out to $\abss{z}=R\rightarrow \infty$ . By doing so, we pick up the residues at each of the zeros $z_k$ of $ze^{-z}-e^{iu-e^{iu}}$, except $z=e^{iu}$ since the latter is a removable pole of the integrand, where

\beqn
\abss{z_1}\leq \abss{z_2} \leq ...\,, z_ke^{-z_k}=e^{iu-e^{iu}},\,\, z_k \neq e^{iu}.
\eeqn
It is straightforward to check that $z=e^{iu}$ is the only solution of modulus one. We calculate the residues and obtain:
\bneqn \label{eq14}
\phi_m(u)=\sum_{k=1}^\infty \frac{1}{z_k^m}{\frac{1-e^{iu}}{1-z_k}} \sim \frac{1}{z_1^m} \frac{1-e^{iu}}{1-z_1},
\eneqn
since $\abss{z_k}>\abss{z_1}$ for $k>1$. 
We now return to equation (\ref{eq15}) and make the substitution $u=\frac{\theta}{\sqrt{m}}$. Then:
\beqn
z_1e^{-z_1}=e^{iu-1-iu+u^2/2}=e^{u^2/2-1},
\eeqn
and thus
\bneqn \label{eq16}
z_1=1+\epsilon,\,\, \epsilon=\pm i\theta/\sqrt{m}.
\eneqn
Placing (\ref{eq16}) into equation (\ref{eq14}), we obtain 
\beqn
\phi_m\parens{\frac{\theta}{\sqrt{m}}}e^{-i\theta\sqrt{m}} \sim e^{-\frac{\theta^2}{2}}, 
\eeqn
which is (\ref{eq15}). This finishes the proof. $\Cox$
\end{proof}

\section{Simple random war} \label{sec3}
We now consider the model in the work of \cite{Kaigh}, in which each army loses a solider independently of the size of $m$ and $n$. Let the success probability be 1/2, independent of $m$ and $n$. We refer to such a model a ``simple random war.'' Without further delay, we note that the results of this section up until equation (\ref{niceone2}) can also be obtained via the framework of a simple symmetric random walk; however, for purposes of consistency, we employ the framework from Section \ref{sec2}. \\
\indent Let $q(m,n)$ denote the probability that army A is reduced to 0 soldiers before Army B is reduced to 0 soldiers.  $q(m,n)$ satisfies the recurrence
\bneqn
q(m,n)=\frac{1}{2}q(m-1,n)+\frac{1}{2}q(m,n-1),
\eneqn
with
\bneqn
q(m,0)=0,\,\, m>0; \,\, q(0,n)=1,\, n>0,
\eneqn
as the boundary conditions.
By induction, we easily have
\bneqn \label{niceone2}
q(m,n)=\sum_{k=0}^{n-1} \frac{1}{2^{m+k}}{m+k-1 \choose k}.
\eneqn
This result is in agreement with \cite{Kaigh} (p. 23, equation (3)). However, the work of \cite{Kaigh} does not consider asymptotics. We now proceed with Theorem \ref{thmscaling}.

\begin{theorem} \label{thmscaling}
As $m \rightarrow \infty$,
\bneqn \label{simple}
q(m, m+x\sqrt{m}) \rightarrow \Phi \parens{\frac{x}{\sqrt{2}}},\,\, -\infty<x<\infty.
\eneqn
\end{theorem}
\begin{proof}
We first define the characteristic functions analogous to (\ref{eq11}) as
\beqn
\Psi_m(u)=\sum_{n=1}^\infty \bracks{q(m,n)-q(m,n-1)}e^{inu}=E(e^{i \eta_m u}), 
\eeqn
where $\eta_m$ is a random variable with
\beqn
P(\eta_m \leq n)=q(m,n).
\eeqn
Summing by Newton's formula,
\beqn
(1-x)^{-m}=\sum_{n=0}^\infty {m+n-1 \choose n}x^n,
\eeqn
we find, letting $x=\exp{iu}$,
\beqn
\Psi_m(u)=e^{iu}(2-e^{iu})^{-m}.
\eeqn
Making a substitution of $u=\theta/\sqrt{m}$, we have,
\bneqn
\Psi_m\parens{\frac{\theta}{\sqrt{m}}}e^{-i\theta \sqrt{m}}=E(e^{i\theta \frac{\parens{\eta_m-m}}{\sqrt{m}}})\rightarrow e^{-\theta^2},\,\, -\infty<\theta<\infty.
\eneqn
Invoking L\'evy's convergence theorem, the result follows.
$\Cox$
\end{proof} 
Table \ref{1} in the Appendix displays values of $p(m,n)$ and $q(m,n)$ for various values of $m$ and $n$.

\section{Scaling limit}\label{sec4} We start with a scaling of the problem. First, assume the random evolution is a pure jump process in continuous time. Let $(X_{t}, Y_{t})_{t\ge 0}$ be the joint process of the 
fortunes at time $t\ge 0$ of both players and $q(x, y)=\frac{x}{x+y}$, $x+y>0$. At times
$t'>0$ separated by independent and identically distributed (i.i.d) exponential waiting times of mean one, the process is updated according to the transition matrix
\begin{align}
\begin{array}{ccccc}\label{tmx}
X_{t'} &=& X_{t'-} \qquad &\text{with probability}&\quad q\\
X_{t'} &=& X_{t'-} -1\qquad &\text{with probability}&\quad 1-q\\
Y_{t'} &=& Y_{t'-} \qquad &\text{with probability}&\quad 1-q\\
Y_{t'} &=& Y_{t'-} -1\qquad &\text{with probability}&\quad q,
\end{array}\,,
\end{align}
where $q=q(X_{t'-}, Y_{t'-})$.
 We shall assume the process is defined on a filtered probability space $(\Omega, {\mathcal F}, P, ({\mathcal F}_{t})_{t\ge 0})$ satisfying the usual conditions. 
Note that
\begin{align}\label{esum}
X_{t}+Y_{t}=X_{0}+Y_{0} - {\mathcal N}_{t}\,,\end{align} where $({\mathcal N}_{t})_{t\ge 0}$ is the Poisson process of intensity one defined by the exponential waiting times. 
 
 The scale at which we define the process is seen as {\em macroscopic} in the sense of (\ref{neceq}) in Theorem \ref{thm1}. More precisely, we introduce a scaling factor $N>0$ such that
\begin{align}
t_{micro}\rightarrow N t_{macro}\,,\quad x_{micro}=N x_{macro}\,,\quad 
y_{micro}=N y_{macro},
\end{align}
such that the quantities $t$, $X_{t}$, $Y_{t}$ from (\ref{tmx}) are 
microscopic, or equivalently, amplified by a factor of $N$. The {\em macroscopic} quantities will survive in the limit as $N\to\infty$ and are of order one.

    At a macroscopic scale, the sum of the two processes decreases in steps of size $1/N$
according to a Poisson process with sped up rate $N$. The time-space scaling is ``Eulerian'' since time and space have the same scaling $x/t =const$, and  is not diffusive
(when $x^{2}/t=const$). It is the difference between the two fortunes of the players that shall scale in a diffusive (equivalently, parabolic) manner.  Theorem \ref{t_det} refers to the process under Eulerian scaling, while in the critical scale of equal initial fortunes, the second order approximation is reflected in Theorem 
\ref{thmd}. For references regarding the approximation of Markov processes with differential equations, the reader may consult the classical text of Ethier and Kurtz (\cite{Ethier}).

Henceforth, let $t$ denote the macroscopic time and we let lower case letters denote the other macroscopic quantities. We also shall suppress the subscript {\em macro}.
Denote
\begin{align}\label{s1}
x^{N}_{t}=N^{-1}X_{N t}\,, \qquad
y^{N}_{t}=N^{-1}Y_{N t}\,,\qquad {\mathcal N}^{N}_{t}=N^{-1}{\mathcal N}_{Nt}.
\end{align}
We assume the initial conditions 
\begin{align}\label{2}
\lim_{N\to \infty}x^{N}_{0}=x_{0}\,,\qquad
\lim_{N\to \infty}y^{N}_{0}=y_{0}\,,
\end{align}
and let
\begin{align}\label{201}
T_{N}=z^{N}_{0}=x^{N}_{0}+y^{N}_{0}\,,\qquad
T=z_{0}=x_{0}+y_{0}\,.
\end{align}
As in the discrete time case, the process ends at time \begin{equation}\label{taun}
\tau_{N}=\inf\{t >0\,|\, x^{N}_{t}\wedge y^{N}_{t}=0\}\,,\quad {\mathcal T}_{N}={\mathcal T}(Nx^{N}_{0}, N y^{N}_{0})=N \tau_{N}\,,
\end{equation}
where the capitalized notation designates the time of ruin before scaling.

Theorem \ref{t_det} shows the criticality of the case $x_{0}-y_{0}=0$ and motivates the necessity of a finer scale for the difference.

\begin{theorem}\label{t_det}
Under (\ref{s1}) and (\ref{2}),
 the processes $(x^{N}_{t})$, $(y^{N}_{t})$, $t\in [0, T)$, converge jointly in probability to the deterministic
solutions of the affine equation
\begin{align}\label{xy}
u'_{t}=\frac{u_{t}}{T-t} -1\,,\qquad u_{t}=\frac{u_{0} T}{T-t} +\frac{1}{2}\frac{(T-t)^{2}-T^{2}}{T-t}\,,\quad 
x_{t}+y_{t}=T-t,
\end{align}
with initial values $u_{0}=x_{0}$, respectively $u_{0}=y_{0}$.
Without  loss of generality, assume that $x_{0}\ge y_{0}$. The time to extinction is
\begin{align}\label{xyt}
\tau&=T-\sqrt{x_{0}^{2}-y_{0}^{2}} \le T \,,\qquad &\text{if}\quad x_{0} \ge y_{0}\,,\\
\notag
\qquad \tau&=T\,, \quad &\text{if}\quad x_{0}=y_{0}\,.
\end{align}
\end{theorem}

\begin{proof}
The proof is given in Section \ref{sec71}. $\Cox$
\end{proof}

\begin{corollary}\label{ctn}
If $x_{0}=y_{0}$, then
$\tau_{N} \convp T$.
\end{corollary}
\begin{proof}
We know that 
the process $(x^{N}_{\cdot}, y^{N}_{\cdot})$ converges in distribution to the deterministic pair of ordinary differential equations $(x_{\cdot}, y_{\cdot})$
given in (\ref{xy}). Fix $t>0$.
For any Borel set $F$, if $P(x_{t}y_{t}\in \partial F)=0$, where $\partial F$ is the boundary of $F$, we have
\begin{align}\label{F}
\lim_{N\to \infty} P(x^{N}_{t}y^{N}_{t}\in F)=P(x_{t}y_{t}\in F).
\end{align}
The product $x_{t}y_{t}$ is deterministic and hence has a delta distribution. If $t <\tau$ from (\ref{xyt}), $x_{t}y_{t}\not=0$ and $0$ is a continuity point of the distribution. Since $\{\tau_{N}\le t\}= \{x^{N}_{t}y^{N}_{t}=0\}$, for $F=\{0\}$
we have shown that
if $t <\tau$, then $\lim_{N\to \infty}P(\tau_{N}\le t)=0$.
When $x_{0}=y_{0}$ and $\tau=T$,  $\tau_{N}\le T_{N}$ almost surely since the process decreases only from the initial value. Since $\lim_{N\to \infty}T_{N}=T$, $\tau_{N} \convd T$. Since $T$ is non-random, the convergence holds in probability. $\Cox$

\end{proof}

\begin{proposition}\label{pt1}
Under the same assumptions as in Theorem \ref{t_det}, and using the same notation as in equation (\ref{eq1}),
if $x_{0}-y_{0}\not=0$,
\begin{align}
\lim_{N\to \infty} p(N x^{N}_{0}, N y^{N}_{0})={\bf 1}_{x_{0}<y_{0}}
\,.
\end{align}
\end{proposition}

\begin{proof}
The result can be proven directly or as a corollary of 
Theorem \ref{thm1}. 
In the continuous time setting,
$p(m, n)$ is defined identically in distribution over the skeleton Markov chain at jump times.
  By noticing that $p(m, \cdot)$ is increasing and $n=Ny^{N}_{0} >> x$ 
as $N\to \infty$
(here $x$ is as in 
Theorem \ref{thm1}),
we see that the limit must be greater than any value $\Phi(x)$, $x\in {\mathbb R}$, and thus it is one. The opposite case is obtained by complementarity. $\Cox$
\end{proof}

It is now clear that $x_{0}-y_{0}=0$ is the critical case. 
We now require some additional assumptions on the 
initial states,
\begin{align}
\label{s2}
z^{N}_{t}=N^{\frac{1}{2}}(x^{N}_{t} - y^{N}_{t})\,,\quad
\lim_{N\to \infty}z^{N}_{0}=z_{0}\,.
\end{align}
The last condition requires us to have \begin{equation}\label{tx}
x_{0}=y_{0}\,\,, \quad T=2x_{0}\,.
\end{equation}

\begin{theorem}\label{thmd}
Under (\ref{s2}) and (\ref{tx}), in addition to the conditions of Theorem \ref{t_det}, the difference process
$(z^{N}_{t})$ converges in distribution, as a process on the 
Skorokhod space of right continuous with left limits path space, to the diffusion
\begin{align}\label{3}
d z_{t}&=\frac{z_{t}}{T-t} dt + d W_{t}\,,\qquad
 0\le t <T\quad
 \text{starting at}\quad z_{0}.
\end{align}
Further, equation (\ref{3}) can be solved explicitly as the Gaussian process,
\begin{align}\label{4}
z_{t}=\left(1-\frac{t}{T}\right)^{-1}
\left[z_{0}+\int_{0}^{t} \Big(1-\frac{s}{T}\Big) d W_{s}\right]
\,.
\end{align}
\end{theorem}

\begin{proof}

The proof is given in Section \ref{sec72}. $\Cox$
\end{proof}
 
\begin{remark}
Note that equation (\ref{3}) is similar to the stochastic differential equation (SDE) satisfied by the Brownian bridge, with the difference that the drift has a positive sign (alternatively, positive direction) relative to the sign of $z_{t}$. This is significant, leading to equation (\ref{4}), which shows that the mean (when $z_{0}\not=0$) and standard deviation of the process $z_{t}$ are of size $(1-\frac{t}{T})^{-1}$. Unlike the Brownian bridge, the mean and standard deviation do not shrink to zero as $t\to T$. In fact, as $t\to T$, $z_{t}\to \infty$ almost surely.
\end{remark}
  
\subsection{$C$-tightness}

We begin with a definition.
\begin{definition}

\indent Let a family of right-continuous with left limits processes $(\eta^{N}_{t})$, defined on $t\in [0, T)$, with values in ${\mathbb R}^{d}$, $d\ge 1$, indexed by $N>0$, be {\em $C$-tight} 
if, for any $T'\in (0, T)$ and $\epsilon>0$, conditions $(i)$ and $(ii)$ below are satisfied. To simplify notation, we let $|\cdot|$ denote the Euclidean norm, irrespective of dimension.

\begin{align}\label{ct1}
(i) \qquad \limsup_{M\to\infty}
\limsup_{N\to \infty}P(\sup_{t\in [0, T']} |\eta^{N}_{t}|>M) =0\,,\\
\label{ct2}
(ii) \qquad \lim_{\delta\to 0}\limsup_{N\to \infty}P(\sup_{t,t' \in [0, T'] \, \,,|t-t'|< \delta} |\eta^{N}_{t}-\eta^{N}_{t'}|] > \epsilon) =0.
\end{align}
\end{definition}
\noindent In the proof of Theorem \ref{t_det} we need briefly to consider $d=2$ for $\eta^{N}_{t}=(x^{N}_{t}, y^{N}_{t})$,
but in the proofs of Theorems \ref{thmd} and \ref{tsc}, we only need that $d=1$.

Condition $(i)$ of the definition states uniform boundedness on any compact time interval and condition $(ii)$ states the uniform modulus of continuity of the family of processes based on the {Arzela-Ascoli Theorem}. $C$-tightness refers to the fact that
 processes, in our case (\ref{s1}), (\ref{s2}), defined as pure jump processes, live on the Skorokhod space $D=D([0, T), {\mathbb R})$, endowed with the $J1$ metric as a Polish space. Conditions $(i)$ and $(ii)$ guarantee that the probability laws $P_{N}$ of 
 $(\eta^{N}_{t})$ are precompact on $D$ and that their limit points are continuous paths in the subspace 
 $C([0, T), {\mathbb R})\subseteq D([0, T), {\mathbb R})$. For further references, the reader may consult \cite{Bill} and \cite{Shiryaev}.

We now let $\phi(x)$ be a real, non-negative function, 
 having at most exponential growth. Possibilities for $\phi(x)$ include all absolute values of polynomials 
 as well as other standard test functions (that is,  $\phi(x)=|x|^{m}$, $m$ even or $\phi(x)=\exp{\lambda x}$, $\lambda>0$). We shall now replace $(i)$ in (\ref{ct1}) with the stronger 
  condition,
  \begin{equation}\label{ct11}
(i') \qquad 
\limsup_{N\to \infty}E[\sup_{t\in [0, T']} \phi(\eta^{N}_{t})] <\infty,
\end{equation}
since condition $(i)$ becomes a consequence of Markov's inequality.   
  
 Let $f\in C_{b}^{1,2}([0, T)\times {\mathbb R}^{d})$ be a test function with bounded derivatives.   
   Recall that the process is sped up by a factor of $N$, appearing in front of the time integrals below. 
   We give (see, for example, \cite{Landim-Kipnis})
  the analogue of It\^{o}'s differential formula for pure jump processes. It states that $({\mathcal M}^{N, f}_{t})$ and 
 $(\overline{{\mathcal M}}^{N, f}_{t})$ are $({\mathcal F}_{t})$ - martingales for $t\in [0, \tau_{N}]$. The jumps $J^{N}_{\pm}(f)$
 and their probabilities given by $q^{N}_{s}$
correspond to the dynamics in (\ref{tmx}), after the scaling 
 (\ref{s1}). They are
 \begin{align}\label{mfg}
   {\mathcal M}^{N, f}_{t}&=
   f(t, x^{N}_{t}, y^{N}_{t})-f(0, x^{N}_{0}, y^{N}_{0})\\
   \notag
   &-\int_{0}^{t}
  N\Big[\partial_{s}f(t, x^{N}_{s}, y^{N}_{s}) + J^{N}_{+}(f) q^{N}_{s} + J^{N}_{-}(f) (1-q^{N}_{s})\Big] \, ds,
  \end{align}
  \text{and}
  \begin{align}\label{mfgj}
   \overline{{\mathcal M}}^{N, f}_{t}= ({\mathcal M}^{N, f}_{t})^{2}
   -\int_{0}^{t}
  N\Big[ (J^{N}_{+}(f))^{2} q^{N}_{s} + (J^{N}_{-}(f))^{2} (1-q^{N}_{s})\Big] \, ds,
  \end{align}
  where
  \begin{align}
\notag
  J^{N}_{+}(f)&= \left(f(t, x^{N}_{s}, y^{N}_{s}-\frac{1}{N}) - f(t, x^{N}_{s}, y^{N}_{s}\right),
  \\
  \notag
   J^{N}_{-}(f)&= \left(f(t, x^{N}_{s}-\frac{1}{N}, y^{N}_{s}) - f(t, x^{N}_{s}, y^{N}_{s}\right)\,,
   \\
   \notag
   q^{N}_{s}&=  \frac{x^{N}_{s}}{x^{N}_{s}+y^{N}_{s}}\,.
  \end{align}

\section{The residual time to extinction in the critical case}
  
Throughout the entirety of this section we shall assume that $x^{N}_{0}+y^{N}_{0}=T_{N}=T$, a slightly
 stronger initial condition than (\ref{201}).
Consistent with the scaling in  (\ref{s2}), we note
that, for a given $N$,
 the process $(z^{N}_{t})$ lives in the angle bounded by 
$|z/(T_{N}-t)|< N^{\frac{1}{2}}$,
 within the small error appearing in (\ref{51}).
Let $\tau_{N}$, defined in (\ref{taun}),
be the 
time, on macroscopic scale, when one of the players reaches zero. Then:
\begin{align}\label{5}
\tau_{N}=\tau_{N}(T, z_{0})&=\inf\{ t >0\, |\, |z^{N}_{t}| = (x^{N}_{t}+y^{N}_{t})N^{\frac{1}{2}}\}\,,\\
\label{51}
x^{N}_{t}+y^{N}_{t}&=T_{N}-{\mathcal N}^{N}_{t}=T_{N}-t +o(N^{-\frac{1}{2}})
\,.
\end{align}

Equations (\ref{s2}) 
and (\ref{tx}) show that the time to ruin is a hitting time of a set depending on the initial state $z^{N}_{0}$ and the total fortune at time zero, given by $T_{N}$. We shall denote $T_{N}$ by $\tau_{N}$ when the other dependence is not essential.
 Before scaling, the time to ruin is
\begin{align}\label{bigTau}
{\mathcal T}_{N}=
{\mathcal T}_{N}(N x^{N}_{0}, N y^{N}_{0})=
N \tau_{N}.
\end{align}
From (\ref{s1})-(\ref{2}),  recall that $X_{0}+Y_{0}=N T_{N}$ and $\lim_{N\to \infty}T_{N}=T>0$. Corollary \ref{ctn}  shows that, in the critical case,
 $\lim_{N\to \infty}(T-\tau_{N})=0$ in probability. 
 A more refined scaling is necessary to estimate the order of magnitude of the
{\em residual time of ruin} $T-\tau_{N}$.   
 Theorem \ref{tsc} will establish the exact power $\beta=\frac{1}{4}$ such that $S_{N}\convd S$, where
\begin{align}\label{Tt}
T - \tau_{N} \sim S_{N} N^{-\beta}\quad \Leftrightarrow
\quad NT-{\mathcal T}_{N} \sim S_{N} N^{1-\beta}\,,\qquad T=2x_{0}\,,
\end{align}
and  $S>0$ has distribution 
$\theta(T, z_{0})$, which will be determined exactly in Corollary \ref{csc}.

To evaluate the distribution of $T-\tau_{N}$,
we need construct a family of martingales {\em for the limiting process $(z_{t})$ from
(\ref{3})}, indexed by a parameter $\rho$ and
adapted to the filtration of the process. This family has the form
$f(T-t,z_{t})$ given by
 \begin{align}\label{6}
f(a, u)=a^{\rho}g(a^{-\frac{1}{2}}u)\,,\qquad a>0\,,\quad u, \rho \in {\mathbb R}\,,
\end{align}
 where 
  \begin{align}\label{7}
g''(u)+3u g'(u) -2\rho g(u)=0\,.
\end{align}
A power series solution to (\ref{6})
 is of the form
  \begin{align}\label{8}
g_{\rho}(u)=\sum_{k\ge 0} a_{k} u^{k}\,,\quad a_{k+2}=\frac{2\rho-3k}{(k+2)(k+1)} a_{k}\,, \qquad k=0, 1. \ldots\,
\end{align}
We note the presence of two independent solutions, one 
 even and one odd, which occurs because the coefficients depend on a two step recurrence. The solution is 
 convergent on the real line and is a linear combination of the odd and even solution, with coefficients $a_{0}$ and $a_{1}$ determined by the initial conditions of the second order ordinary differential equation (\ref{7}).
 
Equation (\ref{7}) can be solved after a substitution
 $g(u)=h\parens{\frac{u^{2}}{2}}$, with $h$ obtained from a Kummer function
 (see \cite{wiki}) in (\ref{k2}). This implies that $g=g_{\rho}$, depending on the exponent $\rho$,
 is even in $u$. 
We have 
 \begin{align}\label{13}
x h''(x) +\parens{\frac{1}{2}+ 3 x}h'(x) - \rho h(x)=0.
\end{align}
 Let $w(z)=M(a, b, z)$ be the solution of the Kummer equation
  \begin{align}\label{k1}
z w''(z) +(b-z) w'(z) - a w(z)&=0\,,\qquad b=\frac{1}{2}\,,\quad a=\frac{1}{2}+\frac{\rho}{3}\,,\\
w(z)&=\sum_{k=0}^{\infty} \frac{a^{(n)}}{b^{(n)} n!} \, z^{n}
=\sum_{k=0}^{\infty} \frac{a^{(n)}}{ (n!)^{2}} \, z^{n}
\,,\notag
\end{align}
 where $a^{(0)}=1$ and for $k\ge 1$ integer
 $a^{(k)}=\Pi_{j=0}^{k-1} (a+j)$ is the {\em rising factorial}. 
 The Kummer transformation
\begin{align}\label{k21}
e^{-z}M(a, b, z)=M(b-a, b, -z),
\end{align}
 gives
\begin{align}\label{k2}
h_{\rho}(x)=e^{-3x} w(3x)=e^{-3x} M\parens{\frac{1}{2}+\frac{\rho}{3}, \frac{1}{2}, 3x}
=M\parens{-\frac{\rho}{3}, \frac{1}{2}, -3x}\,.
\end{align}
It is remarkable that these functions are generalized Laguerre polynomials for $-\frac{\rho}{3}\in {\mathbb Z}$, which will become more transparent after we identify the moments of the limiting random variable $S$. As $z\to \infty$,
 \begin{align}\label{Kinf}
M(a, b, z)\sim \Gamma(b)\left( \frac{e^{z} z^{a-b} }{\Gamma(a)} +  \frac{(-z)^{-a} }{\Gamma(b-a)}\right)\,.
\end{align}
Given (\ref{Kinf}) and the power series form of 
the Kummer function at $z=0$, we have that, for any $\rho>0$, there exist $0<c_{1}(\rho)<c_{2}(\rho)$ such that
 \begin{align}\label{hinf}
 c_{1}(\rho) |z|^{\frac{\rho}{3}}\le h_{\rho}(z)\le c_{2}(\rho)|z|^{\frac{\rho}{3}},\end{align}
for all $z$. 
Working with the Kummer function's integral representation or with the {\em contiguous relations}, there exists a positive $c_{3}(\rho)$ such that
 \begin{align}\label{hinf1}
  |h^{(k)}_{\rho}(z)|\le c_{3}(\rho)(|z|\vee 1)^{\frac{\rho}{3}-k}\,,\quad k\ge 0\,.
 \end{align}
The above being established, we are ready to present the main result of this section. It identifies 
the distribution 
$\theta(T, z_{0})$ on $(0, \infty)$ of the scaled asymptotic residual time, defined in equation (\ref{Tt}).

\begin{theorem}\label{tsc}
If, in addition to the conditions of Theorem \ref{thmd}, the initial 
value $T_{N}=T$, then
the scaled residual time $S_{N}=N^{\frac{1}{4}}(T-\tau_{N})$ converges in distribution as $N\to \infty$ to a positive random variable $S$. Its distribution $\theta(T, z_{0})$
has the exact generating function 
\begin{align}\label{tsc1}
E[S^{q}]=\parens{\frac{2}{3}}^{\frac{q}{4}}  \frac{\Gamma(\frac{1}{2}+\frac{q}{4})}{\Gamma(\frac{1}{2})} T^{\frac{3q}{4}}h_{\frac{3q}{4}}\parens{\frac{z_{0}^{2}}{2T}}\,,\quad q>0\,.
\end{align}
The above fully determines the distribution on $(0, \infty)$.
\end{theorem}
\begin{proof}
The result is proven step by step in Section \ref{s:tsc}. 
Regarding the identification of $S$, formula (\ref{tsc1}) is established in subsection \ref{sbs:1} and the positive definiteness of $S$ is proven in subsection \ref{sbs:2}. 
$\Cox$
\end{proof}

\medskip\noindent
\begin{remark}
We note that (\ref{tsc1}) is obtained from (\ref{k2}) by setting
$q=\frac{4\rho}{3}$.
\end{remark}

In fact, the fourth power of $S$, modulo a constant, has a classical distribution.

\begin{corollary}\label{csc}
The distribution 
of the adjusted asymptotic residual time 
$R=3T^{-3} S^{4}$
is the {\em non-central chi squared-distribution} with 
$k=1$ degrees of freedom and 
non-centrality parameter $\frac{3z_{0}^{2}}{T}$.
In the case $z_{0}=0$, we have $h_{\rho}(0)=1$ and 
$R$ is $\chi^{2}_{1}$, i.e. $R$ is the square of a standard normal.
\end{corollary}

\begin{proof}
Set $p$, $q$ such that
$m=\frac{q}{4}=\frac{\rho}{3}$, $m\ge 0$ an integer.
Using (\ref{k21}) and (\ref{k2}) 
we see that  
\begin{equation}\label{chi01}
h_{3m}(w)=M\parens{-m, \frac{1}{2}, -3w}= 
\frac{m! \Gamma(\frac{1}{2})}{\Gamma\parens{\frac{1}{2}+m}}
L_{m}^{\parens{-\frac{1}{2}}}(-3w)\,,\qquad w=\frac{z_{0}^{2}}{2T}\,,\end{equation}
where $L_{m}^{\parens{-\frac{1}{2}}}(x)$, $x\in {\mathbb R}$ is the $m$-th generalized Laguerre polynomial of type $\alpha=-\frac{1}{2}$, the latter defined as
$$L_{m}^{(\alpha)}(x)=
\frac{x^{-\alpha}e^{x}}{n!}\frac{d^{m}}{d x^{m}}(x^{m+\alpha} e^{-x})=
{m +\alpha \choose m}M(-m, \alpha+1,x)\,.$$ The generating function $\sum_{m=0}^{\infty}\lambda^{m}
L_{m}^{(\alpha)}(x)$, when 
$x=-3w$, $\alpha=-\frac{1}{2}$,
combining with (\ref{tsc1}),
is then equal to the moment generating function of a random variable with
moments
given as \begin{equation}\label{chi02}
m!L_{m}^{\parens{-\frac{1}{2}}}(-3w)
= \frac{\Gamma\parens{\frac{1}{2}+m}}{\Gamma\parens{\frac{1}{2}}}h_{3m}(w)=M\parens{-m, \frac{1}{2}, -3w}\,.\end{equation}
This is equal (cf. \cite{wiki-r}) to
\begin{align}\label{chisq}
M(\lambda)=\sum_{m=0}^{\infty} \lambda^{m}L_{m}^{\parens{-\frac{1}{2}}}(-3w)
=\frac{e^{\frac{
3w \lambda}{1-\lambda}}}{(1-\lambda)^{\frac{1}{2}}}\,,
\end{align}
which is equivalent to saying it is the distribution of $\frac{1}{2}R'$, where $R' \sim \chi^{2}_{1}(6w)$, the
non-centered $\chi^{2}$ distribution with two degrees of freedom and non-centrality parameter 
$6w=\frac{3z_{0}^{2}}{T}$. Alternatively, it is the square of a normal with nonzero mean with variance one. 

The proof is straightforward. It immediately follows by combining (\ref{chi01}), (\ref{chi02}), and (\ref{chisq}) with (\ref{tsc1}). $\Cox$
\end{proof}

\begin{corollary}\label{tcont}
 The residual time after extinction before scaling 
$NT-{\mathcal T}_{N}=O(N^{\frac{1}{4}})$
in the precise sense
 \begin{align}\label{ct10}
N^{-\frac{3}{4}} [TN - {\mathcal T}_{N}(N x_{0} , N x_{0} - N^{\frac{1}{2}} z_{0})] \convd S\,.
\end{align}
\end{corollary}

\begin{proof}
This is an immediate consequence of Theorem \ref{tsc} and the scaling in (\ref{s1}). $\Cox$
\end{proof}

 \section{Proof of Theorem \ref{tsc}}\label{s:tsc}
 
The goal of this section is to prove Theorem \ref{tsc}, the second key result of our paper. To do so, we first state and prove Propositions \ref{prepl}, \ref{p:num}, and \ref{pab}. Throughout, let the scaled residual times be
\begin{align}\label{sN}
\hat{S}_{N}=N^{\beta}(T-{\mathcal N}^{N}_{\tau_{N}})\,,\quad S_{N}=N^{\beta}(T-\tau_{N})\,,\qquad \beta=\frac{1}{4}\,.
\end{align}

\subsection{Preliminary bounds}

\begin{proposition}\label{prepl}
The following inequality holds:
$$
E[|\hat{S}_{N}-S_{N}|^{2}] \le TN^{-\frac{1}{2}}.
$$
Further, if one of the families of random variables
$(\hat{S}_{N})_{N>0}$,  $(S_{N})_{N>0}$ is tight with limit $S$ in distribution, then so is the other.
\end{proposition}
 
 \begin{proof}
Doob's maximal inequality applied to the square-integrable martingale ${\mathcal N}^{N}_{t}-t$ (the compensated Poisson process scaled by $N$) satisfies
\begin{align}\label{doob}
 E[\sup_{s\in [0, T]} |{\mathcal N}^{N}_{\tau_{N}} - \tau_{N}|^{2}]
\le T N^{-1}\,.
\end{align}
We obtain the desired result by taking the expected value of $|\hat{S}_{N}-S_{N}|^{2}$, noting that the supremum dominates the value at $t=\tau_{N}$, and multiplying with $N^{2\beta}$. The exponent of $N$ is $2\times \frac{1}{4}-1=-\frac{1}{2}$. The second statement of the proposition is immediate. $\Cox$
 \end{proof}

 \begin{proposition}\label{p:num}
 Let $a,b > 0$ with $\frac{b}{a} > \frac{\ln 2}{\ln \frac{3}{2}}$.
 For any $x\ge y\ge 1$
 \begin{align}\label{inab}
\Big[(x+y-1)^{a} &(x-y+1)^{b}-(x+y)^{a}(x-y)^{b}\Big]x\\
 \notag
&+ \Big[(x+y-1)^{a} (x-y-1)^{b}-(x+y)^{a}(x-y)^{b}\Big]y \ge 0\,.
 \end{align}
  \end{proposition}
  
  \begin{proof}
  Divide both sides by $(x+y-1)^{a}(x-y)^{b}(x+y)^{-1}$. The inequality becomes
 $$ \frac{x}{x+y}\Big(1+\frac{1}{x-y}\Big)^{b}+ \frac{y}{x+y}\Big(1-\frac{1}{x-y}\Big)^{b} \ge  \Big(1+\frac{1}{x+y-1}\Big)^{a}\,.$$
 Let $f(u)=u^{b}$, a convex function for $u\ge 0$.
 Jensen's inequality implies that it is sufficient to show
 $$ \Big(1+\frac{1}{x+y}\Big)^{b}\ge
  \Big(1+\frac{1}{x+y-1}\Big)^{a}\,.
  $$ For $z=\frac{1}{x+y}$ the above expression is equivalent to
  $$(1+z)^{b}(1-z)^{a}\ge 1\,, \quad 0<z<\frac{1}{2}\,,\qquad \frac{b}{a}
   > \frac{\ln 2}{\ln \frac{3}{2}}\,.$$
  The logarithm $g(z)=b\ln(1+z) + a\ln (1-z)$ is well defined since $0<z<\frac{1}{2}$ has $g(0)=0$ and in all cases, $g'(0)>0$. 
Since the critical point is $z_{0}=\frac{b-a}{b+a}$,
it is the case that $z_{0}> \frac{1}{2}$ when $b/a>3$ and thus $g'(0)>0$ on $z\in (0, \frac{1}{2}]$. This shows the above inequality holds. If $b/a \le 3 $, the minimum will be achieved at either endpoint. We only have to verify the inequality at $z=\frac{1}{2}$, where the function takes the value $b\ln \frac{3}{2} - a \ln 2$, which is non-negative by assumption.
$\Cox$
  \end{proof}

 \begin{proposition}\label{pab}
Let $a,b > 0$ as in Proposition \ref{p:num}.
Then $H^{N}_{t}(a,b)=(x^{N}_{t}+y^{N}_{t})^{a}(z^{N}_{t})^{b}$, where $z^{b}=(z^{2})^{\frac{b}{2}}$, is a ${\mathcal F}_{t}$-submartingale. Using Doob's maximal inequality for $p>1$,
\begin{align}\label{smg1}
E[(\sup_{s\in [0, \tau_{N}]}H^{N}_{s}(a, b))^{p}] \le N^{pc}\parens{\frac{p}{p-1}}^{p}E[\hat{S}_{N}^{(a+b)p}]\,,\quad c=\frac{b-a}{4}=\frac{b}{2}-\frac{(a+b)}{4}\,.
\end{align}
 \end{proposition} 
 
 \begin{proof}
We base the proof on  (\ref{inab}). After dividing by $N$ everywhere, we notice that the $ds$ term in
(\ref{mfg}) applied to $H^{N}_{t}(a,b)$ gives the left hand side of the expression in (\ref{inab}). If the integrand is non-negative, $H^{N}_{t}(a,b)$ is a sub-martingale.
Without loss of generality, assume $x\ge y$.  By construction, we stop at $\tau_{N}$ when the minimum is zero, but we know it is impossible that both $x$ and $y$ equal zero at the same time. 
With these remarks, 
to verify that the conditions of (\ref{inab}) are met, we also note that the expression under the integral is calculated only {\em before $\tau_{N}$}, i.e. $Nx^{N}_{t}\ge N y^{N}_{t}\ge 1$
for $t< \tau_{N}$. 
We then simply apply
(\ref{inab}), proving the sub-martingale claim. We shall have an inequality between the values of the submartingale at time $0$ and $\tau_{N}$. 
The stopping time is uniformly bounded since $\tau_{N}\le T$, for all $N$.
We employ Doob's maximal inequality in $L^{p}$ norm. 
Taking into account (\ref{5}), the inequality is easily verified. We conclude by verifying exponent of $N$. At $\tau_{N}$ a factor of $N^{\frac{1}{2}}$ enters the factor with exponent $b$ only. To compensate the power $N^{\beta(a+b)}$, where $\beta=\frac{1}{4}$, we must subtract, giving the formula for $c$ in (\ref{smg1}). $\Cox$

 \end{proof}

 \subsection{Plan of the proof of Theorem \ref{tsc}}
 We shall apply (\ref{ito00}) to the function $f(x^{N}_{t}+y^{N}_{t}, z^{N}_{t})$ from (\ref{6}). Ideally we would like to work with
$f(T-t, z^{N}_{t})$, but this would add errors that would hinder 
the proof. Notice that $z^{N}_{t}=\sqrt{N}(x^{N}_{t}-y^{N}_{t})$, so that the application of the It\^{o} formula is direct in the form
given in the setup for the process $(x^{N}_{t},y^{N}_{t})$.
On the other hand, $x^{N}_{t}+y^{N}_{t}=T_{N}-{\mathcal N}^{N}_{\tau_{N}}$, and Proposition \ref{prepl}
will be necessary to 
replace $T_{N}-{\mathcal N}^{N}_{\tau_{N}}$ with $T-\tau_{N}$.

The family of functions $f$, for all $\rho\in {\mathbb R}$, constructed 
to satisfy (\ref{7}) and (\ref{8}), will nullify 
the time integral in It\^{o}'s formula for $(z_{t})$, making it a martingale for the limiting diffusion. 

Of course, we are using the formula for $(z^{N}_{t})$, and an error term is expected. After showing this error term is vanishing as $N\to \infty$, uniformly in $t\in [0, \tau_{N}]$, we apply an optional stopping argument to evaluate the martingale at both ends $t=0$ and $t=\tau_{N}$. 

\subsection{The differential formula}
The function $g_{\rho}$ (\ref{8}) that defines $f$ relates to
the Kummer function $h=h_{\rho}$ associated to $\rho$ by $g(u)=h(\frac{u^{2}}{2})$, present in (\ref{k2}).
It\^{o}'s formula (\ref{mfg})-(\ref{mfgj}) in this case
will give the terms
\begin{align}\notag
\left[f\Big(x^{N}_{s}+y^{N}_{s} -\frac{1}{N}, \frac{z^{N}_{s}+\frac{1}{\sqrt{N}}}{\sqrt{x^{N}_{s}+y^{N}_{s} -\frac{1}{N}}}\Big)-
f\Big(x^{N}_{s}+y^{N}_{s}, \frac{z^{N}_{s}}{\sqrt{x^{N}_{s}+y^{N}_{s} }}\Big)\right]\times \frac{x^{N}_{s}}{x^{N}_{s}+y^{N}_{s}}
\end{align}
\begin{align}\notag
+\, \left[f\Big(x^{N}_{s}+y^{N}_{s} -\frac{1}{N}, \frac{z^{N}_{s}-\frac{1}{\sqrt{N}}}{\sqrt{x^{N}_{s}+y^{N}_{s} -\frac{1}{N}}}\Big)-
f\Big(x^{N}_{s}+y^{N}_{s}, \frac{z^{N}_{s}}{\sqrt{x^{N}_{s}+y^{N}_{s} }}\Big)\right]\times \frac{y^{N}_{s}}{x^{N}_{s}+y^{N}_{s}}\,.
\end{align}

All terms present in the Taylor formula of order two - we include the second derivative to match (\ref{ito7}) - including the error terms in Lagrange form, 
depend on the first three derivatives of $f$ with respect to $z$
and the second derivative with respect to $t$. They
are of the form
$a^{\rho-i}u^{i'} h_{\rho}^{(k)}(\frac{u^{2}}{2a})$,$\,0\le k\le 3$ and $0\le i, i'\le 4$, with $a=x^{N}_{s}+y^{N}_{s}$ and $u=z^{N}_{s}$.

At $z\to \infty$, the Kummer functions satisfy (c.f. \cite{wiki})
(\ref{Kinf}) and (\ref{hinf})
and, based on (\ref{10}) and the choices of $a$, $b$ in (\ref{k2}),
the function $h^{(k)}_{\rho}(z)$, and by consequence $z\to f(T-t, z)$
are of polynomial order $\frac{\rho}{3}-k$
at $z\to \infty$. 
After inspecting carefully the derivatives, we see that
in the asymptotic formula, modulo constants independent of $N$, $a$, $u$, but possibly not $\rho$,
\begin{equation}\label{er01}
|\partial_{aa}f(a, u)| \le c_{21}(\rho) a^{\frac{2\rho}{3}-2} \, u^{\frac{2\rho}{3}}\,,\qquad 
\partial^{3}_{z}f(a, u) \le c_{22}(\rho) a^{\frac{2\rho}{3}} \, 
u^{\frac{2\rho}{3}-3}\,.
\end{equation}
Additionally, the term $\partial_{aa}f$ has a factor of $N\times (\frac{1}{N})^{2}$ and the term $\partial^{(3)}_{z}f$ has a factor of $N\times (\frac{1}{\sqrt{N}})^{3}$.
This shows that (\ref{ito00}) is a martingale plus an error term ${\mathcal E}^{N}(f)={\mathcal E}^{N}_{1}(f)+{\mathcal E}^{N}_{2}(f)$.

\subsection{The error terms from the space derivatives and time derivatives}
We start with the error term issued from the time derivatives present in (\ref{er01}),

\begin{equation}\label{enn1}
{\mathcal E}^{N}_{1}(f)
\le
N^{-1}
T C_{21}(\rho)\sup_{s\in [0,\tau_{N}]}H^{N}_{s}(\frac{\rho}{3}-j, \frac{4\rho}{3}-j')\,,\quad j=2\,, j'=0\,.
\end{equation}
We continue with the error term issued from the space derivatives present in (\ref{er01}),
\begin{equation}\label{enn2}
{\mathcal E}^{N}_{2}(f)
\le
N^{-\frac{1}{2}}
T C_{22}(\rho)\sup_{s\in [0,\tau_{N}]}H^{N}_{s}(\frac{\rho}{3}-j, \frac{4\rho}{3}-j')\,,\quad j=0\,, j'=3
\,.\end{equation}
For $\rho>\frac{9}{4}$, let $$p=\frac{\frac{2\rho}{3}+\frac{2\rho}{3}}{\frac{2\rho}{3}+\frac{2\rho}{3}-j-j'}>1\,.$$
Above, we have written the sum of two terms $2/3$ because one comes from the power of $a=T-\tau_{N}$, which is $\rho-\rho/3=2\rho/3$, while the power of $u=z$ is $2\times \rho/3$. The exponent $\rho/3$ is the asymptotic order of magnitude of the Kummer functions $h_{\rho}$.

\subsection{Bounds based on Proposition \ref{pab}}
The error comes from from the derivatives in time and space of $f$ from Taylor's formula. 
Then, with $C_{21}(\rho)$, $C_{22}(\rho)$ constants incorporating usual terms in the Taylor series of order up to three
as well as $c_{21}(\rho)$, $c_{22}(\rho)$, 
with the factor $T$ coming from the time integral, we 
calculate, with $p=p_{1}$ corresponding to the second derivative in time, i.e. $j=2$, $j'=0$,
\begin{align}\label{cp1}
(E[({\mathcal E}^{N}_{1}(f))^{p}])^{\frac{1}{p}}
&\le
N^{-1}
T C_{21}(\rho)
(E[\sup_{s\in [0,\tau_{N}]}(H^{N}_{s}(\frac{2\rho}{3}-j, \frac{2\rho}{3}-j'))^{p}] )^{\frac{1}{p}}\\
\notag
&\le N^{-1+c_{1}}
T C_{21}(\rho)(E[\hat{S}_{N}^{\frac{4\rho}{3}}])^{\frac{1}{p}},
\end{align}
with $c_{1}$ calculated for $a=\frac{2\rho}{3}-j$, $b=\frac{2\rho}{3}-j'$
and $j=2$, $j'=0$, yielding
\begin{align}
\notag
c_{1} -1=  \frac{j-j'}{4}-1=-\frac{1}{2}
\,.\end{align}
Similarly, we calculate, with $p=p_{2}$ corresponding to the third derivative in space, i.e. $j=0$, $j'=3$,
\begin{align}\label{cp2}
(E[({\mathcal E}^{N}_{2}(f))^{p}])^{\frac{1}{p}}
&\le
N^{-\frac{1}{2}}
T C_{22}(\rho)
(E[\sup_{s\in [0,\tau_{N}]}(H^{N}_{s}(\frac{\rho}{3}-j, \frac{4\rho}{3}-j'))^{p}] )^{\frac{1}{p}}\\
\notag
&\le N^{-\frac{1}{2}+c_{2}}
T C_{22}(\rho)(E[\hat{S}_{N}^{\frac{5\rho}{3}}])^{\frac{1}{p}},
\end{align}
with $c_{2}$ calculated for $a=\frac{2\rho}{3}-j$, $b=\frac{2\rho}{3}-j'$, and 
giving for $j=0$, $j'=3$,
\begin{align}
\notag
c_{2} -\frac{1}{2}=
\frac{j-j'}{4}-\frac{1}{2}=-\frac{5}{4}
\,.\end{align}

\subsection{Equation based on the Optional Stopping Theorem}

For fixed $N$, $\tau_{N}< T$ almost surely because not
both $x^{N}_{t}$, $y^{N}_{t}$ reach zero at the same time and $T_{N}-\tau_{N}=x^{N}_{\tau_{N}}+y^{N}_{\tau_{N}}$.
We are ready to apply the Optional Stopping Theorem
to equate the initial value with the final value at $t=\tau_{N}$ and obtain
\begin{align}\label{90}
E[(T-{\mathcal N}^{N}_{\tau_{N}})^{\rho} g_{\rho}(
\frac{z^{N}_{\tau_{N}}}{\sqrt{T-{\mathcal N}^{N}_{\tau_{N}}}})]= T^{\rho}g_{\rho}(\frac{z_{0}}{\sqrt{T}})+E[{\mathcal E}_{N}(f)]\,.
\end{align}
As we mentioned in
(\ref{5}), at $\tau_{N}$ we have the identity 
$$|z^{N}_{\tau_{N}}|
=\sqrt{N}(T_{N}-{\mathcal N}^{N}_{\tau_{N}})\approx \sqrt{N}(T_{N}-\tau_{N}),$$
obtained from algebraic manipulation of (\ref{s2}). The replacement 
of the Poisson process with the value $\tau_{N}$ will be done based on Doob's maximal inequality, given that $\tau_{N}\le T$ uniformly in $N>0$, as in Proposition \ref{prepl}.
We have used that $g(u)$ is even, being a function of $u^{2}$, having equal values at $z=\pm (T-t)$. These correspond to player $X$ (with $+$) respectively $Y$ (with $-$) being the winners.
Relation (\ref{90}) can be written in terms of
the Kummer function $h=h_{\rho}$ associated to $\rho$ by $g(u)=h(\frac{u^{2}}{2})$ (\ref{k2}) as follows:
\begin{align}\label{10}
E[(T-\tau_{N})^{\rho} h_{\rho}(N(T-\tau_{N}))]= T^{\rho}h_{\rho}(\frac{z_{0}^{2}}{2T})+E[{\mathcal E}_{N}(f)]\,.
\end{align}

\subsection{Tightness of $(\hat{S}_{N})_{N>0}$ and $(S_{N})_{N>0}$}

Denote
$v_{N}=E[\hat{S}_{N}^{\frac{4\rho}{3}}]$. Notice that
(\ref{hinf}) gives a lower bound 
$$c_{1}(\rho) v_{N} \le E[(T-{\mathcal N}^{N}_{\tau_{N}})^{\rho} g_{\rho}(
\frac{z^{N}_{\tau_{N}}}{\sqrt{T-{\mathcal N}^{N}_{\tau_{N}}}})],$$
and that we have the upper bound as
$$
 E[{\mathcal E}_{N}(f)]\le 
N^{-\frac{1}{2}}
T C_{21}(\rho)v_{N}^{\frac{1}{p_{1}}}
+
N^{-\frac{5}{4}}
T C_{22}(\rho)v_{N}^{\frac{1}{p_{2}}}\,.
$$ 
This implies that, with the obvious meaning of the constants,
$$c_{1}(\rho)
v_{N} \le T^{\rho}g_{\rho}(\frac{z_{0}}{\sqrt{T}})+ \bar{c}_{0}N^{-\frac{1}{2}} + 
\bar{c}_{1}N^{-\frac{1}{2}} (v_{N}+1),$$
showing that 
$$\limsup_{N\to \infty}v_{N}=\limsup_{N\to \infty}
E[\hat{S}_{N}^{\frac{4\rho}{3}}]\le  c_{1}(\rho)^{-1}T^{\rho}g_{\rho}(\frac{z_{0}}{\sqrt{T}})< \infty\,.$$
This proves that $(\hat{S}_{N})$ is a tight family of non-negative random variables.
Let $S$ be a limit point in distribution sense. We know from
Proposition \ref{prepl} that $(S_{N})$ is also tight.

\subsection
{Distribution of the positive part $S_{+}$ of the limiting distribution}\label{sbs:1}

We proceed to prove that any limit point $S$ has distribution completely determined by (\ref{tsc1}), and as a consequence it is unique in law. We conclude that $S_{N}\convd S$.
To simplify notation, we keep the same notation $\hat{S}_{N}$ for the subsequence converging to $S$. 
Fix a sufficiently small $\epsilon>0$. 
First we determine the moments of $S\vee \epsilon$.

Returning to (\ref{10}), we now know that the error term vanishes as $N\to \infty$ and
\begin{align}
E[(T-{\mathcal N}^{N}_{\tau_{N}})^{\rho} h_{\rho}(N(T-{\mathcal N}^{N}_{\tau_{N}}))]
=N^{-\frac{\rho}{4}}E[\hat{S}_{N}^{\rho} h_{\rho}(N^{\frac{3}{4}}\hat{S}_{N})]
\\
= N^{-\frac{\rho}{4}}E[(\hat{S}_{N}\vee \epsilon)^{\rho} h_{\rho}(N^{\frac{3}{4}}(\hat{S}_{N}\vee \epsilon))] + C(N,\epsilon),
\end{align}
where 
\begin{align}
C(N,\epsilon)\le c_{2}(\rho) \epsilon^{\frac{4\rho}{3}}\,.
\end{align}

At this point, 
our goal is to determine the distribution of the limit $S$. 
As we see, the last relation involves only the deterministic functions $h_{\rho}$ and $\hat{S}_{N}$.
We shall appeal to 
the Skorokhod representation theorem (see \cite{Bill1}).
 There exists a probability space supporting all $(S_{N})$, $S$, such that all preserve the same distributions and $S_{N}\to S$ almost surely.

The random variables $\hat{S}_{N}\vee \epsilon\ge \epsilon$ and $N^{\frac{3}{4}}(\hat{S}_{N}\vee \epsilon) \to \infty$ pointwise. Due to this fact and (\ref{Kinf}) we know exactly the limit:
 $$\lim_{N\to \infty}
 \frac{h_{\rho}(N^{\frac{3}{4}}(\hat{S}_{N}\vee \epsilon))}{(N^{\frac{3}{4}}(\hat{S}_{N}\vee \epsilon))^{\frac{\rho}{3}}}=
\Big(
(\frac{2}{3})^{\frac{\rho}{3}}
\frac{\Gamma(\frac{1}{2}+\frac{\rho}{3})}{\Gamma(\frac{1}{2})}\Big)^{-1}\quad a.s.$$
Then
$$
\lim_{N\to \infty}N^{-\frac{\rho}{4}}E[(\hat{S}_{N}\vee \epsilon)^{\rho} h_{\rho}(N^{\frac{3}{4}}(\hat{S}_{N}\vee \epsilon))]
=
\Big(
(\frac{2}{3})^{\frac{\rho}{3}}
\frac{\Gamma(\frac{1}{2}+\frac{\rho}{3})}{\Gamma(\frac{1}{2})}\Big)^{-1}
E[(S\vee \epsilon)^{\frac{4\rho}{3}}],
$$
which proves by dominated convergence that
if $S_{+}=S {\bf 1}_{S>0}$ then
\begin{align}\label{s+}
E[S_{+}^{\frac{4\rho}{3}} ]=\lim_{\epsilon\downarrow 0}
E[(S\vee \epsilon)^{\frac{4\rho}{3}}]
=
(\frac{2}{3})^{\frac{\rho}{3}}
\frac{\Gamma(\frac{1}{2}+\frac{\rho}{3})}{\Gamma(\frac{1}{2})}
T^{\rho}h_{\rho}(\frac{z_{0}^{2}}{2T})\,,\qquad \forall \rho>\frac{9}{4} \,,
\end{align}
which is essentially (\ref{tsc1}).

\subsection{Uniqueness of the moments and $S_{+}=S$}
\label{sbs:2}
The next step is to identify that a random variable with moments 
given by (\ref{s+}) 
is uniquely determined if the moments are known on 
an interval $(3, \infty)$. This is true 
due to the Mellin transform theorem (see \cite{Widder},  Theorem 6a, p 243).
Finally we want to show that a distribution with those moments
has exactly the moments on the right hand side of (\ref{s+}), {\em for all $\rho>0$}.
This is proven in the derivation of (\ref{chisq}).

The final step is to show that such a distribution 
does not concentrate at zero. This is true because now we can use all moments down to $q=0$ for $S_{+}$. For $q=0$ we obtain the zero-th moment, equal to the total mass of $S_{+}$. This is equal to one, which shows that
$S=S_{+}$ or $P(S=0)=0$.

\subsection{Special case $z_{0}=0$}
In case $z_{0}=0$
we obtain that $S^{2}\sim c_{0} \, |Z|$, $Z\sim N(0, 1)$, $c_{0}$ constant determined in Corollary \ref{csc}.
This concludes the proof of
Theorem \ref{tsc}. \hfill $\Box$

\section{Proofs of Theorems \ref{t_det} and \ref{thmd}} \label{sec7}
\subsection{Proof of Theorem \ref{t_det}} \label{sec71}

 The processes $(x^{N}_{t})$, $(y^{N}_{t})$ are bounded in the interval $[0, T]$ by construction, so condition 
  $(i)$ in (\ref{ct1}) immediately holds. The test function can be assumed simply continuous in all derivatives since the time and state spaces are compact. Let $f(t, x, y)=x$, while $f(t, x, y)=y$ is analogous. 
  Equations (\ref{mfg})-(\ref{mfgj}) show that there exists
  $m^{N}(t', t)= O(\frac{1}{N})$ uniformly in $t, t'$ such that
  $$x^{N}_{t}-x^{N}_{t'}=\int_{t'}^{t}(1-q^{N}_{s}) ds +
 m^{N}(t', t)(t-t').
 $$
This shows that condition (ii) in (\ref{ct2}) holds. Conditions $(i)-(ii)$ guarantee that any limit point of the tight sequence is actually continuous. 
 
  For a general test function $f(t, x, y)=f(t,x)$, let \begin{equation}\label{opa}
  A_{t}f(t,x)=\frac{x}{T-t}\partial_{x}f(t,x)-1,
  \end{equation}
  corresponding to  (\ref{3}). 
 Given $(\eta_{t})$ a path in the Skorohkod space and $t\in [0, T']$, $0<T'<T$, we obtain that the functional $\Psi$,
 $$\eta \to \Psi(\eta)
  =\sup_{t\in [0, T']}\Big|
  f(t, \eta_{t})- f(0, \eta_{0})-\int_{0}^{t}
  \partial_{s}f(s, \eta_{s}) + A_{s}f(s, x_{s})
  \, ds\Big|,$$
 is continuous and bounded.
  Let $x_{\cdot}$ be a limit point of the tight sequence $(x^{N}_{\cdot})_{N>0}$. 
   The quadratic variation is easily verified
  to be of order $1/N$ according to (\ref{mfgj}).
  To explain the operator $A_{t}$ in (\ref{opa}), note that
  the first order term in the Taylor series at $x=x^{N}_{s}$ with change $-\frac{1}{N}$,
  is $$\partial_{x}f(s, x^{N}_{s}) (1-q^{N}_{s})= 
  \partial_{x}f(s, x^{N}_{s}) \,\Big(\frac{T-{\mathcal N}^{N}_{s}-x^{N}_{s}}{T-{\mathcal N}^{N}_{s}} \Big)
  =A_{s} f(s, x^{N}_{s})+ e_{N},$$
  with $E[e_{N}^{2}]\to 0$ uniformly in time.
  The error term is controlled by (\ref{doob}).

  From the continuity theorem,
  $E[\Psi(x)]=0$.
  This reasoning 
  is repeated in more detail in the proof of Theorem \ref{thmd} in Section \ref{sec72}.
  It follows that $x_{\cdot}$ is a 
  path (possibly random) satisfying, almost surely,
  $$ f(t, x_{t})- f(0, x_{0})-\int_{0}^{t}
  \partial_{s}f(s, x_{s}) + 
  L_{s}f(s, x_{s})\, ds=0\,.$$
  Since $x_{\cdot}$ is continuous, the equality must be satisfied for all $t$. 
  There is only one possible continuous solution to this ordinary differential equation posed in integral (weak) form. One can verify the existence and uniqueness of the strong solution of 
  (\ref{xy}) by standard results for differential equations
  with continuous, Lipschitz coefficients on $[0, T']$, $0<T'<T$. The solution as well as the extinction time are elementary.
  \hfill
  $\Box$

\subsection{Proof of Theorem \ref{thmd}} \label{sec72}

  First, we recall that due to (\ref{esum}) and (\ref{s2}), the joint process $(z^{N}_{t}, {\mathcal N}^{N}_{t})$ is a Markovian pure jump process.
We start by writing the differential equation
 for the process 
$(z^{N}_{t})$ from (\ref{s2}), following 
(\ref{mfg})-(\ref{mfgj}).
  Let $f\in C^{1,2}([0, \infty)\times {\mathbb R}, {\mathbb R})$ be a test function. Note that the processes $(x^{N}_{t})$, $(y^{N}_{t})$, and $(z^{N}_{t})$
are bounded by $N^{\frac{1}{2}}$ (the first two are simply bounded of order one, uniformly in $N$)
almost surely, and so 
the condition on the finite expected value of all martingales in the It\^{o} formula given below is automatically satisfied {\em as long as $N$ remains fixed}. In the limit we shall work with 
the diffusion (\ref{3}) which is not uniformly bounded, even on a compact time interval $[0, T']$, $0< T'< T$.
  
 Appearing in front of the integrals below, the factor $N$ marks the scaling of the time variable.
  The corresponding It\^{o} formula (\ref{mfg}) is
  \begin{align}\label{ito00}
  f(t, z^{N}_{t}) &- f(0, z^{N}_{0}) =N\, \int_{0}^{t} (I) +(II)+(III) \, ds + {\mathcal M}^{N,f}_{t},
  \end{align}
  where
  \begin{align}
  (I) &=
   \partial_{s}f(s, z^{N}_{s}), \nonumber \\
  (II)&= q^{N}_{s}
  \left[f(s, z^{N}_{s}+\frac{1}{\sqrt{N}})-f(s, z^{N}_{s})\right], \nonumber \\
  (III)&=
  (1-q^{N}_{s}) 
 \left[ f(s, z^{N}_{s}-\frac{1}{\sqrt{N}})-f(s, z^{N}_{s})\right],\nonumber \\
 q^{N}_{s}&=\frac{x^{N}_{s}}{x^{N}_{s}+y^{N}_{s}}=\frac{x^{N}_{s}}{ T_{N} - {\mathcal N}^{N}_{s}}\,, \nonumber
  \end{align}
with the last term ${\mathcal M}^{N, f}_{t}$ from (\ref{ito00})
  satisfying that the two processes
   \begin{align}\label{sqj}
 {\mathcal M}^{N,f}_{t}\qquad\text{as well as}\quad
& ({\mathcal M}^{N,f}_{t})^{2}-N\, \int_{0}^{t}
 (J1) + (J2) \, ds\,,\quad 0\le t\le \tau_{N}\\
 (J1)&=q^{N}_{s}
  \left[f(s, z^{N}_{s}+\frac{1}{\sqrt{N}})-f(s, z^{N}_{s})\right]^{2}, \\
  (J2)&=
  (1-q^{N}_{s}) 
 \left[ f(s, z^{N}_{s}-\frac{1}{\sqrt{N}})-f(s, z^{N}_{s})\right]^{2},
  \end{align}
  are $({\mathcal F}_{t})$- martingales. The following five steps, labeled as Steps 1-5, complete the proof.

 \medskip
 {\it Step 1.}
  First, we shall prove an estimate for the factor $(x^{N}_{t}+y^{N}_{t})^{-1}$
appearing in the denominator of the probability $q$ from (\ref{tmx}).  
Fix $T'\in (0, T)$. 

Let $\phi(x)$ be a positive function. Note that the process runs until $\tau_{N}\le T_{N}=x^{N}_{0}+y^{N}_{0}$, making the denominator $x^{N}_{t}+y^{N}_{t}=T_{N}-N^{-1}{\mathcal N}_{N t}$ positive at all times. We proceed to replace $t$ by $t \wedge \tau_{N}$ and obtain:
  
\begin{align}
&E\left[\phi((  x^{N}_{t\wedge \tau_{N}}+y^{N}_{t\wedge \tau_{N}})^{-1})\right]=E\left[\phi(\frac{N}{NT_{N}-{\mathcal N}_{Nt}})\right] \nonumber\\
 &\le \sum_{m=0}^{NT_{N}-1} 
 \phi( \frac{N}{NT_{N}-m}) P({\mathcal N}_{Nt}=m) \nonumber\\
 &\le 
 \phi( \frac{N}{NT_{N}-N(\frac{t+T}{2})}) P({\mathcal N}_{Nt}\le N (\frac{t+T}{2}))
 + \phi(N) P({\mathcal N}_{Nt}> N (\frac{t+T}{2}))\nonumber\\
 &\le 
 \phi(\frac{1}{N^{-1}T_{N}-(\frac{t+T}{2})})+
 \Phi(N)e^{-I(\frac{t+T}{2}) N}\nonumber \,,
  \end{align}
  where $I(\cdot)$ is the large deviations rate function for the Poisson process ${\mathcal N}_{Nt}$.
  For any exponential $\phi$ dominated by the large deviations rate function $I(\cdot)$ for the linear Poisson process of intensity one, the inequality leads to
  an upper bound 
  \begin{align}\label{nn4}
  E\left[\phi((  x^{N}_{t\wedge \tau_{N}}+y^{N}_{t\wedge \tau_{N}})^{-1})\right] \le C(\phi, t),
  \end{align}
  where the constant $C(\phi, t)$ does not depend on $N>0$
  but $\lim_{t\to T}C(\phi, t)=+\infty$. 
  Setting $t=T'$ we obtain (\ref{nn4}) on $[0, T']$.
  
  \medskip
  {\em Step 2.}
  Let $f(t, z)=z$ in It\^{o} formula. It follows that
  \begin{align}\label{ito1}
  z^{N}_{t}=z^{N}_{0} + N \int_{0}^{t} \frac{1}{\sqrt{N}} \frac{x^{N}_{s}-y^{N}_{s}}{T_{N}-{\mathcal N}^{N}_{s}}\, ds +
  {\mathcal M}^{N,z}_{t},
 \end{align}
  and, by (\ref{s2}), the integrand is equal to $z^{N}_{s}/(T_{N}-{\mathcal N}^{N}_{s})$. 
  It is straightforward that (\ref{ito1}) gives directly 
  (\ref{3}) since the martingale part converges, in distribution, to a Brownian motion, as an application of 
the invariance principle to a scaled compensated Poisson process.

However, this does not close the argument.
It is easy to square the equality and obtain after employing Cauchy-Schwartz, H{\"o}lder's inequality for the time integral, as well as the bound $t\le T$, that
    \begin{align}\label{ito2}
 \frac{1}{3} [z^{N}_{t}]^{2}\le 
  [z^{N}_{0}]^{2} +  T \int_{0}^{t}  [\frac{z^{N}_{s}}{T_{N}-{\mathcal N}^{N}_{s}}]^{2}\, ds +
  [{\mathcal M}^{N,z}_{t}]^{2}.
 \end{align}
  
  We recall that we only need to work with a compact subinterval $[0, T']\subseteq [0, T)$. 
  Combining this with (\ref{nn4}) for $\phi(x)=x^{2}$ and denoting the constant by $C(2, T')$ since we shall only concentrate on the time interval $[0, T']$, we obtain
    \begin{align}\label{ito3}
 \frac{1}{3} [z^{N}_{t}]^{2}\le 
  [z^{N}_{0}]^{2} +  TC(2,T')  \int_{0}^{t}  [z^{N}_{s}]^{2}\, ds +
  [{\mathcal M}^{N,z}_{t}]^{2}.
 \end{align}
  After inspecting the quadratic variation part of (\ref{sqj})  we see that we can bound the last square term in (\ref{ito3}) by $C_{2}T$. We do so by Doob's maximal inequality, obtaining
    \begin{align}\label{ito4}
 \frac{1}{3} [z^{N}_{t}]^{2}\le 
  [z^{N}_{0}]^{2} +  T C(2,T')  \int_{0}^{t}
   \sup_{0\le s'\le s}[z^{N}_{s}]^{2}\, ds +
  C_{2}T.
 \end{align}
 We now take the supremum over time 
  on the right hand side, apply Doob's maximal inequality to the martingale part, and observe that Gr{\"o}nwall's lemma
  gives the bound
    \begin{align}\label{ito5}
 E [\sup_{0\le t\le T'} \phi(z^{N}_{t})]\le 
  3([z^{N}_{0}]^{2} + C_{2}T)
  e^{ T C(2,T') T'}:=C_{3}(T')\,,\qquad \phi(z)=z^{2}\,,
 \end{align}
 which is (\ref{ct11}). 
 
\begin{remark}
The dependence on $T'<T$ is essential in the bound since $C(2,T')$ blows up as $T'\uparrow T$.
\end{remark}
 
 \medskip
 {\em Step 3.}
 We shall obtain a similar bound as (\ref{ito5})
 for exponential functions. We start with $\phi(z)=e^{\lambda z}$, $\lambda>0$. The difficult term in It\^{o}'s formula is the time integrand
 $$N e^{\lambda z^{N}_{s}}\left[(e^{\frac{\lambda}{\sqrt{N}}} -1) \frac{x^{N}_{s}}{T_{N}-{\mathcal N}^{N}_{s}}
 +
 (e^{-\frac{\lambda}{\sqrt{N}}} -1) \frac{y^{N}_{s}}{T_{N}-{\mathcal N}^{N}_{s}}
 \right],$$
 where the prefactor $N$ is from the time scaling. 
Given that $|e^{h}-1-h|\le c_{4} h^{2} e^{|h|}$ 
and
considering again the general bound (\ref{nn4}) from Step 1,
we 
see that an upper bound  for
the above term reduces to a a constant depending on $T'$ times
\begin{equation}\label{dtbd1}
 \lambda e^{\lambda z^{N}_{s}} [z^{N}_{s} + c_{4}\lambda^{2}]\,. \end{equation}
We have used Cauchy-Schwarz's followed by Doob's maximal inequality.
The martingale part is of a smaller order of magnitude. Putting together the bound obtained in (\ref{dtbd1}) with the bound on the quadratic variation, and employing analogous notation to that of Step 2, together with  (\ref{ito5}), we obtain 
  \begin{align}\label{ito6}
 E [\sup_{0\le t\le T'} e^{2\lambda z^{N}_{t})}]\le 
  3(\phi(z^{N}_{0}) + C_{5}(T))
  e^{ C(4,T') T'}:=C_{4}(T')\,. 
 \end{align}

\medskip
{\em Step 4.}
To verify the modulus of continuity (ii) in (\ref{ct2}),
we subtract (\ref{ito1}) at two times $s,t$ and observe that
under the time integral we have, after simplification with the scaling constants, a term bounded above by
$$E\left[\sup_{s, t \in [0, T']\,,\, |t-s|<\delta}
\left|
\int_{s}^{t} \frac{z^{N}_{s'}}{T_{N}-{\mathcal N}^{N}_{s'}} \, ds'\right|
\right]\le \delta E\left[(
\frac{1}{T_{N}-{\mathcal N}^{N}_{T'}})
\sup_{s' \in [0, T']}|z^{N}_{s'}|\right]\,,
$$
where we used the observation that the function $s'\to
(T_{N}-{\mathcal N}^{N}_{s'})^{-1}$ is nondecreasing. 
Cauchy-Schwarz's inequality, with (\ref{nn4}) for $\phi(x)=|x|^{2}$ and the bound 
(\ref{ito5}) prove (ii) for the time integral term. The martingale part is straightforward due to Doob's maximal inequality applied to the second moment. We have thus proven that $(z^{N}_{t})$ is tight.

\medskip
{\em Step 5.}
Keeping $T'\in (0, T)$ we notice that the limit process (\ref{3}) is uniquely defined by its martingale problem. This is true for less regular coefficients but here we have smooth, bounded, Lipschitz coefficients. Thus, to identify the limit, we proceed in the standard way. 
First, we assume $(z_{t})$ is a continuous path process 
found as a possible limit point of the tight sequence $(z^{N}_{t})$. 
Second, we shall show that for any 
 $f\in C^{1,2}_{c}([0, \infty)\times {\mathbb R}, {\mathbb R})$
as in (\ref{ito00}), only with compact support for convenience (this is a determining class for the martingale problem),
\begin{align}\label{ito7}
f(t, z_{t})-f(0, z_{0})-
\int_{0}^{t}\left(
\partial_{s}f(s, z_{s}) + 
\frac{z_{s}}{T-s} \partial_{z}f(s, z_{s})+\frac{1}{2} 
\partial_{zz}f(s, z_{s})\right)
\, ds,
\end{align}
 is a $({\mathcal F}_{t})$ martingale. 
 Comparing (\ref{ito7}) with (\ref{ito00}), we see that
 when we replace $z_{t}$ with $z^{N}_{t}$ in (\ref{ito7})
 we obtain a martingale ${\mathcal M}^{N, f}_{t}$ 
 plus an error term ${\mathcal E}_{N}$.

The error term needs to go to zero uniformly in $N$ in 
 $L^{1}$ norm. However, the error term is less complicated here. 
The function $f$ has compact support and thus satisfies a uniform bound in $N$ and $t\in [0, T']$. Simply applying Taylor's formula in Lagrange form with degree two - and error of degree three -
 $$\Big|
 f(s, z+\epsilon) -f(s, z)- \epsilon \partial_{z}f(s, z)
 -\frac{\epsilon^{2}}{2}\partial_{zz}f(s, z)\Big| \le C(f) \epsilon^{3}\,,\quad \epsilon=\frac{1}{\sqrt{N}}\,.
 $$
 Notice that the terms containing $\epsilon$ are exactly the 
 second order operator corresponding to diffusion (\ref{3}). 
 Terms (II) and (III) in (\ref{ito00}) give the error
 $${\mathcal E}_{N} \le N \int_{0}^{t} C(f) N^{-\frac{3}{2}} ds \le N^{-\frac{1}{2}} C(f)T'\,.$$
The rest of the proof is technical but standard. 
 If $\Phi(\eta(\cdot))$ is a continuous bounded functional of a path $\eta\in D([0, T'], {\mathbb R})$ and $(z^{N}_{\cdot})\convd (z_{\cdot})$, the {\em continuity theorem} implies that
 $$\lim_{N\to \infty}E[\Phi(z^{N}_{\cdot})]=
E[\Phi(z_{\cdot})]\,.$$
In our case, we fix $s,t$, $0\le s\le t\le T'$ and a bounded function $\Gamma(\eta_{\cdot})$ measurable with respect to ${\mathcal F}_{s}$. Then by writing
$$
\Phi(\eta_{\cdot}):=[f(t, \eta_{t}) - f(s, \eta_{s}) - \int_{s}^{t}
L_{s'} f(s', \eta_{s'}) \, ds']\Gamma(\eta_{\cdot}),
$$
we can verify that the functional $\Phi$ is continuous, bounded in the $J1$ topology 
on the Skorokhod space. It has expected value equal to the error term exactly, and when we pass to the limit we obtain that the expected value with respect to the law of the limiting process vanishes. It follows that
$(z_{t})$ satisfies the martingale problem (\ref{ito7}), and by uniqueness, it is the solution to (\ref{3}). 
%

%
%

\medskip
\noindent
{\bf Acknowledgments.}
We would like to thank Larry Shepp and Robert W. Chen for 
introducing the problem to us. We also wish to thank Larry Shepp for his kind assistance with some parts of Section \ref{sec2} of the manuscript and Min Kang for helpful discussions. We thank Quan Zhou for his excellent assistance with numerics. Finally, we are extremely grateful to an anonymous reviewer whose invaluable suggestions greatly improved the quality of this work.

\newpage
\noindent
{\bf Appendix}\\
We display the values of $p(m,n)$ and $q(m,n)$ in Table \ref{1} below. 
\newcolumntype{C}[1]{>{\centering\let\newline\\\arraybackslash\hspace{0pt}}m{#1}}
\begin{table}[h!]
\centering
\noindent
\begin{tabular}{C{1.5cm} C{1.5cm} C{1.5cm} C{2cm} C{2cm} C{2cm} C{2cm} } 
\toprule
\multicolumn{3}{c}{Initial Capital} & \multicolumn{2}{c}{Win probability} & \multicolumn{2}{c}{Win probability} \\ 
\cmidrule(lr){1-3}
\cmidrule(lr){4-5}
\cmidrule(lr){6-7}
$m+n$ &  $m$  &  $n$ &  $p(m, n)$ & $1 - p(m,n)$ & $q(m,n)$ & $1 - q(m, n)$  \\ 
\cmidrule(lr){1-7}
  &  8 & 12  & 0.939 & 0.061 & 0.820 & 0.180 \\
20 & 9 & 11 & 0.779 & 0.221 & 0.676 & 0.324 \\
 & 10 & 10 & 0.5 & 0.5 & 0.5 & 0.5 \\
\cmidrule(lr){1-7}
  &  45 & 55  & 0.958 & 0.042 & 0.843 & 0.157 \\
100 & 48 & 52 & 0.755 & 0.245 & 0.656 & 0.344 \\
 & 50 & 50 & 0.5 & 0.5 & 0.5 & 0.5 \\
\cmidrule(lr){1-7}
  &  90 & 110  & 0.993 & 0.007 & 0.922 & 0.078 \\
200 & 95 & 105 & 0.890 & 0.110 & 0.761 & 0.239 \\
 & 100 & 100 & 0.5 & 0.5 & 0.5 & 0.5 \\
\cmidrule(lr){1-7}
  &  480 & 520  & 0.986 & 0.014 & 0.897 & 0.103 \\
1000 & 490 & 510 & 0.863 & 0.137 & 0.737 & 0.263 \\
 & 500 & 500 & 0.5 & 0.5 & 0.5 & 0.5 \\
\cmidrule(lr){1-7}
  &  960 & 1040  & 0.999 & 0.001 & 0.963 & 0.037 \\
2000 & 980 & 1020 & 0.939 & 0.061 & 0.815 & 0.185 \\
 & 1000 & 1000 & 0.5 & 0.5 & 0.5 & 0.5 \\
\bottomrule
\end{tabular}
\caption{Values of $p(m,n)$ and $q(m,n)$.}\label{1} 
\end{table}

\end{document}